\documentclass[11pt]{amsart}
\usepackage{amsmath}
\usepackage{amssymb}
\usepackage{latexsym}
\usepackage{graphicx}

\newtheorem{theorem}{Theorem}%[section]
\newtheorem{corollary}{Corollary}%[section]
\newtheorem{lemma}{Lemma}%[section]
\theoremstyle{remark}
\newtheorem{remark}{Remark}%[section]
\begin{document}

\title[coefficient bounds]{\large Bounds on the coefficients of certain analytic and univalent functions}

\author[K. O. Babalola]{K. O. BABALOLA}

\begin{abstract}
For the real number $\alpha>1$, we use a technique due to Nehari
and Netanyahu and an application of certain integral iteration of
Caratheodory functions to find the best-possible upper bounds on
the coefficients of functions of the class $T_n^\alpha (\beta)$
introduced in \cite{TOO} by Opoola.
\end{abstract}

%\begin{amssubject}
%30C45.
%\end{amssubject}

%\begin{keyword}
%Coefficient bounds, analytic and univalent functions.
%\end{keyword}

\maketitle

\section{Introduction}
Let $A$ denote the class of functions:
\begin{equation}
f(z)=z+a_2z^2+...\, \label{1}
\end{equation}
which are analytic in the unit disk $E=\{z\in \mathbb{C}\colon
|z|<1\}$. In \cite{TOO}, Opoola introduced the subclass
$T_n^\alpha (\beta)$ consisting of functions $f \in A$ which
satisfy:
\begin{equation}
Re \frac{D^nf(z)^\alpha}{\alpha^nz^\alpha}>\beta\, \label{2}
\end{equation}
\\where $\alpha>0$ is real, $0\le\beta<1$, $D^n(n\in N_0=\{0, 1, 2,
...\})$ is the Salagean derivative operator defined as:
$D^nf(z)=D(D^{n-1}f(z))=z[D^{n-1}f(z)]^{\prime}$ with
$D^0f(z)=f(z)$ and powers in ~(\ref{2}) meaning principal
determinations only. The geometric condition ~(\ref{2}) slightly
modifies the one given originally in \cite{TOO} (see \cite{BO}).

The object of the present work is the extension of some earlier
results regarding the bounds on the coefficients, $a_k$, of
functions belonging to the class $T_n^\alpha (\beta)$. Babalola
and Opoola have begun to solve this problem in \cite{BO2}. They
determined sharp bounds on $a_k$ for some $\alpha>0$ and gave some
rough estimate of the general coefficient bounds using the
logarithmic coefficient approach. The sharp bounds were stated as
follows:
\begin{theorem}[\cite{BO2}]
Let $f\in T_n^\alpha (\beta)$. Define
\[
B_m=\frac{2^m(1-\beta)^m\alpha^{m(n-1)}\prod_{j=0}^{m-1}(1-j\alpha)}{m!}
\]
\[
A_k(n,\alpha,\beta)=\sum_{m=1}^{k-1}B_mQ_{k-1}^{(m)}
\]
\[
A_{k-1}(n,\alpha,\beta)=\sum_{m=1}^{k-2}B_mQ_{k-1}^{(m)}
\]
where for each $m=1, 2, ..., Q_{k-1}^{(m)}$ is defined by the
power series
\[
\left(\sum_{k=1}^{\infty}\frac{z^k}{(\alpha+k)^n}\right)^m=Q_{m}^{(m)}z^m+Q_{m+1}^{(m)}z^{m+1}+Q_{m+2}^{(m)}z^{m+2}+...
\]
Also let
\[
\Omega_1=\{\alpha|0<\alpha<(k-2)^{-1}, k=2, 3, ...\},
\]
\[
\Omega_2=\{\alpha|(k-2)^{-1}\le\alpha\le (k-3)^{-1}, k=4, 6,...\},
\]
\[
\Omega_3=\{\alpha|(k-2)^{-1}\le\alpha<(k-3)^{-1}, k=3, 5, ...\}
\]
Then
\[
|a_k|\le \left\{
\begin{array}{ll}
A_k&\mbox{if $\alpha\in \Omega_1\cup \Omega_2$},\\
A_{k-1}&\mbox{if $\alpha\in \Omega_3$.}
\end{array}
\right.
\]
The inequalities are sharp. Equalities are attained for $f(z)$
satisfying
\[
\frac{D^nf(z)^\alpha}{\alpha^nz^\alpha}= \left\{
\begin{array}{ll}
\frac{1+(1-2\beta)z}{1-z}&\mbox{if $\alpha\in \Omega_1\cup \Omega_2$},\\
\frac{1+(1-2\beta)z^2}{1-z^2}&\mbox{if $\alpha\in \Omega_3$.}
\end{array}
\right.
\]
\end{theorem}
The rough estimate was also given as:
\begin{theorem}[\cite{BO2}]
Let $f\in T_n^\alpha (\beta), n\in N$. Suppose
\[
\frac{f(z)^\alpha}{z^\alpha}=\sum_{k=0}^{\infty}A_{k+1}(\alpha)z^k,
A_1(\alpha)=1.
\]
Then
\[
|A_{k+1}(\alpha)|<\exp\left\{0.624\alpha^2+(2\alpha^2-\frac{1}{2})\sum_{j=1}^k\frac{1}{j}\right\}.
\]
\end{theorem}
From Theorem 1, which is best-possible, it is obvious that the
problem has only been solved completely for $a_2$ and $a_3$ for
all values of the index $\alpha$. For $k\ge 4$, the problem has
remained open for all $\alpha>(k-3)^{-1}$. In this article we
proceed with the proof of the sharp bounds on the coefficients of
functions in the class $T_n^\alpha (\beta)$, $\alpha>1$. Our
result is the following:
\begin{theorem}
Let $f\in T_n^\alpha (\beta)$. If $\alpha>1$, we have the sharp
inequalities
\begin{equation}
|a_k|\le
\frac{2(1-\beta)\alpha^{n-1}}{(\alpha+k-1)^n},\;\;\;k=2,\;3,\;4,...\,
\label{3}
\end{equation}
Equalities are attained for $f(z)$ satisfying
\[
\frac{D^nf(z)^\alpha}{\alpha^nz^\alpha}=
\frac{1+(1-2\beta)z^{k-1}}{1-z^{k-1}}.
\]
\end{theorem}
For $k=2,\;3$, the above result is contained in Theorem 1 above.
The proof here is however new. Also Singh \cite{RS} proved the
same result for $k=2,\;3,\;4$ for the particular case $n=1$ and
$\beta=0$. In our proof we combine a method of classical analysis
due to Nehari and Netanyau \cite{NN} (also used by Singh
\cite{RS}) with an application of certain integral iteration of
the Caratheodory functions \cite{BO}. That is presented in Section
3.
In the next section we state and prove some preliminary lemmas.

\section{Preliminary Lemmas}
Let $P$ be the class of functions
\begin{equation}
p(z)=1+b_1z+b_2z^2+...\, \label{4}
\end{equation}
which are analytic in $E$ and have positive real part. In
\cite{BO} the following integral iteration of each $p\in P$ was
identified:
\medskip

 {\sc Definition 1(\cite{BO}).} Let $p\in P$ and $\alpha>0$ be real. The
{\em nth} iterated integral transform of $p(z), z\in E$ is defined
as
\[
p_n(z)=\frac{\alpha}{z^\alpha}\int_0^zt^{\alpha-1}p_{n-1}(t)dt,\;\;\;
n\geq 1
\]
with $p_0(z)=p(z)$. The family of the {\em nth} integral iteration
of $p\in P$ was denoted by $P_n$. Functions $p_n(z)$ in $P_n$ have
series expansion:
\begin{equation}
p_n(z)=1+\sum_{k=1}^\infty\frac{\alpha^n}{(\alpha+k)^n}b_kz^k\,.
\label{5}
\end{equation}
Furthermore if $Re\;p_n(z)>\beta,\;(0\le\beta<1)$, we denote by
$P_n(\beta)$ the family of such functions given by:
\[
p_n(z)=1+(1-\beta)\sum_{k=1}^\infty\frac{\alpha^n}{(\alpha+k)^n}b_kz^k.
\]

In the proof of our main result we require the following lemmas:

\begin{lemma}[\cite{NN}]
If $p(z)=1+\sum_{k=1}^\infty b_kz^k$ and $q(z)=1+\sum_{k=1}^\infty
c_kz^k$ belong to $P$, then $r(z)=1+\tfrac{1}{2}\sum_{k=1}^\infty
b_kc_kz^k$ also belongs to $P$.
\end{lemma}

\begin{lemma}[\cite{NN}]
Let $h(z)=1+\sum_{k=1}^\infty d_kz^k$ and
$1+G(z)=1+\sum_{k=1}^\infty b^{\prime}_kz^k$ be functions in $P$.
Set
\begin{equation}
\gamma_m=\frac{1}{2^m}\left[1+\frac{1}{2}\sum_{\mu=1}^m\binom{m}{\mu}d_{\mu}\right],\;\;\;\gamma_0=1\,.
\label{6}
\end{equation}
If $A_k$ is defined by
\[
\sum_{m=1}^{\infty}(-1)^{m+1}\gamma_{m-1}G_1^m(z)=\sum_{k=1}^{\infty}A_kz^k,
\]
then
\[
|A_k|\le 2,\;\;\;k=1,\;2,\;...
\]
\end{lemma}
If, in the proof of the above lemma (as contained in \cite{NN}),we
define $h_n(z)$ as the {\em nth} iterated integral transform of
$h_0(z)=h(z)$ we immediately obtain the following corollary.
\begin{corollary}
Let $h_n(z)$ be the $nth$ integral iteration of
$h_0(z)=1+\sum_{k=1}^\infty d_kz^k$ with $Re\;h_n(z)>\beta$, and
$1+G(z)=1+\sum_{k=1}^\infty b^{\prime}_kz^k$ be functions in $P$.
Define $\gamma_m$ as in ~(\ref{6}) and
\begin{equation}
\eta_m=\frac{(1-\beta)\alpha^n}{(\alpha+m)^n}\gamma_m,\;\;\;\eta_0=1-\beta\,.
\label{7}
\end{equation}
If $A_k$ is defined by
\begin{equation}
\sum_{m=1}^{\infty}(-1)^{m+1}\eta_{m-1}G_1^m(z)=\sum_{k=1}^{\infty}A_kz^k\,,
\label{8}
\end{equation}
then
\begin{equation}
|A_k|\le
\frac{2(1-\beta)\alpha^n}{(\alpha+k)^n},\;\;\;k=1,\;2,\;...\,
\label{9}
\end{equation}
\end{corollary}
\begin{lemma}[\cite{BO2}]
Let $G(z)=\sum_{k=0}^\infty c_kz^k$ be a power series. Then the
mth integer product of G(z) is $G^m(z)=\sum_{k=0}^\infty
c_k^{(m)}z^k$ where $c_k^{(1)}=c_k$ and $c_k^{(m)}=\sum_{j=0}^k
c_jc_{k-j}^{(m-1)}$, $m\ge 2$.
\end{lemma}

We now turn to the proof of the main result.

\section{Proof of the Main result}
Let $f\in T_n^\alpha (\beta)$. Then there exists an analytic
functions $p_n\in P_n$ such that
\begin{equation}
f(z)^\alpha=z^\alpha[\beta+(1-\beta)p_n(z)]\, \label{10}
\end{equation}
where $p_n(z)$ given by ~(\ref{5}) is the {\em nth} iterated
integral transform of an analytic function $p\in P$ defined by
~(\ref{4}) (see Lemma 4.2 of \cite{BO}). The first part of the
proof involves obtaining expression for the coefficients, $a_k$,
of $f(z)$ in terms of the coefficients, $b_k$, of the function
$p\in P$. This is contained in \cite{BO2} and adapted here for
completeness. We have, using ~(\ref{5}) in ~(\ref{10}),
\begin{equation}
f(z)=z\left(1+(1-\beta)\sum_{k=1}^\infty\frac{\alpha^n}{(\alpha+k)^n}b_kz^k\right)^{\frac{1}{\alpha}}\,.
\label{11}
\end{equation}
Equation ~(\ref{11}) expands binomially as
\begin{multline}
\frac{f(z)}{z}=1+\frac{(1-\beta)\alpha^n}{\alpha}\sum_{k=1}^\infty\frac{b_kz^k}{(\alpha+k)^n}+\frac{(1-\beta)^2\alpha^{2n}(1-\alpha)}{2!\alpha^2}\left(\sum_{k=1}^\infty\frac{b_kz^k}{(\alpha+k)^n}\right)^2+\\
...+\frac{(1-\beta)^m\alpha^{mn}\prod_{j=0}^{m-1}{(1-j\alpha)}}{m!\alpha^m}\left(\sum_{k=1}^\infty\frac{b_kz^k}{(\alpha+k)^n}\right)^m+...\,.
\label{12}
\end{multline}
Using Lemma 3 in ~(\ref{12}) we get
\[
\frac{f(z)}{z}=1+\sum_{k=1}^\infty
\widetilde{B}_1C_k^{(1)}z^k+...+\sum_{k=1}^\infty
\widetilde{B}_mC_k^{(m)}z^k+...
\]
where
\[
\widetilde{B}_m=\frac{(1-\beta)^m\alpha^{m(n-1)}\prod_{j=0}^{m-1}{(1-j\alpha)}}{m!}
\]
and $C_k^{(m)}$, $m=1, 2, ...$; $k=m, m+1, ...$ is defined by
\[
\left(\sum_{k=1}^\infty\frac{b_kz^k}{(\alpha+k)^n}\right)^m=\sum_{k=1}^\infty
C_k^{(m)}z^k
\]
having the general form
\begin{equation}
C_k^{(m)}=\sum_{j=1}^k
C_j\prod_{l=1}^m\frac{b_l^{\rho_l}}{(\alpha+l)^{n{\rho_l}}}\,
\label{13}
\end{equation}
for some nonnegative constants $C_j$, $j=1,2,...k$ and indices
$\rho_l$, $l=1, 2, ...,m$ taking values in the set
$M=\{0,1,2,...,m\}$ such that $\rho_l+\rho_2+...+\rho_m=m$ (see
page 10 of \cite{BO2}). From ~(\ref{12}) we write
\begin{equation}
f(z)=z+\sum_{k=2}^\infty \widetilde{A}_k^{(m)}z^k\, \label{14}
\end{equation}
where
\[
\widetilde{A}_k^{(m)}=\sum_{m=1}^{k-1}\widetilde{B}_mC_{k-1}^{(m)},\;\;\;k=2,\;3,\;...
\]
Comparing coefficients in ~(\ref{1}) and ~(\ref{14}), we see that
\[
a_k=\widetilde{A}_k^{(m)}
\]
which gives
\begin{equation}
a_k=\sum_{m=1}^{k-1}\frac{(1-\beta)^m\alpha^{m(n-1)}\prod_{j=0}^{m-1}{(1-j\alpha)}}{m!}\left(\sum_{j=1}^{k-1}
C_j\prod_{l=1}^m\frac{b_l^{\rho_l}}{(\alpha+l)^{n{\rho_l}}}\right)\,.
\label{15}
\end{equation}
Now we compute the leading coefficients, $A_v$, in the expression
~(\ref{8}). From ~(\ref{8}) we have
\begin{equation}
\begin{split}
\sum_{m=1}^\infty
(-1)^{m+1}\eta_{m-1}G_1^m(z)&=G_1(z)-\eta_1G_1^2(z)+...\\
&=\sum_{v=1}^\infty A_vz^v
\end{split}\, \label{16}
\end{equation}
with
\[
G_1(z)=\sum_{v=1}^\infty b^{\prime}_vz^v.
\]
Using Lemma 3 again we have
\begin{equation}
G_1^m(z)=\sum_{v=m}^\infty C_v^{(m)}z^v,\;\;\;m=1,\;2,\;...\,
\label{17}
\end{equation}
where
\begin{equation}
C_v^{(m)}=\sum_{j=1}^v C_j\prod_{l=1}^mb{'}_l^{\rho_l}\,
\label{18}
\end{equation}
with $C_j$ and indices $\rho_l$ as already defined for
~(\ref{13}). Using ~(\ref{17}) and ~(\ref{18}) in ~(\ref{16}) we
get
\[
\begin{split}
\sum_{m=1}^\infty
(-1)^{m+1}\eta_{m-1}G_1^m(z)&=\sum_{v=1}^\infty\left(\sum_{m=1}^v
(-1)^{m+1}\eta_{m-1}C_v^{(m)}\right)z^v\\
&=\sum_{v=1}^\infty\left(\sum_{m=1}^v
(-1)^{m+1}\eta_{m-1}\left(\sum_{j=1}^v
C_j\prod_{l=1}^mb{'}_l^{\rho_l}\right)\right)z^v\\
&=\sum_{v=1}^\infty A_vz^v
\end{split}
\]
so that
\[
A_v=\sum_{m=1}^v (-1)^{m+1}\eta_{m-1}\left(\sum_{j=1}^v
C_j\prod_{l=1}^mb{'}_l^{\rho_l}\right)z^v.
\]
By corollary 1, these coefficients, $A_v$, satisfy the inequality
~(\ref{9}) if $1+G(z)=1+b^{\prime}_1z+b^{\prime}_2z^2+...$ is a
function of the class $P$, and by Lemma 1 we may set
$b^{\prime}_l=\tfrac{1}{2}b_lc_l$ where $p(z)=1+b_1z+b_2z^2+...$
is the function ~(\ref{4}) and $H(z)=1+c_1z+c_2z^2+...$ is an
arbitrary function in $P$. Then
\begin{equation}
|A_v|=\left|\sum_{m=1}^v
(-1)^{m+1}\frac{\eta_{m-1}}{2^m}\left(\sum_{j=1}^v
C_j\prod_{l=1}^mb_l^{\rho_l}c_l^{\rho_l}\right)\right|\le\frac{2(1-\beta)\alpha^n}{(\alpha+v)^n}\,.
\label{19}
\end{equation}
Using ~(\ref{7}) together with the fact that $\alpha>0$, we can
write ~(\ref{19}) equivalently as
\begin{equation}
|A_v|=\left|\sum_{m=1}^v
(-1)^{m+1}\frac{(1-\beta)\alpha^n}{(\alpha+m-1)^n}
\frac{\gamma_{m-1}}{2^m\alpha}\left(\sum_{j=1}^v
C_j\prod_{l=1}^mb_l^{\rho_l}c_l^{\rho_l}\right)\right|\le\frac{2(1-\beta)\alpha^{n-1}}{(\alpha+v)^n}\,.
\label{20}
\end{equation}
Since for each $m=1,2,...$ and any $\alpha>0$,
\[
\prod_{l=1}^m\frac{\alpha^{n{\rho_l}}}{(\alpha+l)^{n{\rho_l}}}\le\prod_{l=1}^m\frac{\alpha^{n{\rho_l}}}{(\alpha+1)^{n{\rho_l}}}=\frac{\alpha^{mn}}{(\alpha+1)^{mn}}\le\frac{\alpha^n}{(\alpha+m-1)^n},
\]
it is evident that for each $v=1, 2, ...$
\begin{multline*}
\sum_{m=1}^v (-1)^{m+1}
\frac{\gamma_{m-1}}{2^m\alpha}\left(\sum_{j=1}^v
C_j\prod_{l=1}^m\frac{(1-\beta)^{\rho_l}\alpha^{n{\rho_l}}}{(\alpha+l)^{n{\rho_l}}}b_l^{\rho_l}c_l^{\rho_l}\right)\le\\
\sum_{m=1}^v (-1)^{m+1}\frac{(1-\beta)\alpha^n}{(\alpha+m-1)^n}
\frac{\gamma_{m-1}}{2^m\alpha}\left(\sum_{j=1}^v
C_j\prod_{l=1}^mb_l^{\rho_l}c_l^{\rho_l}\right).
\end{multline*}
Therefore by ~(\ref{20}),we have
\[
\left|\sum_{m=1}^v (-1)^{m+1}
\frac{\gamma_{m-1}}{2^m\alpha}\left(\sum_{j=1}^v
C_j\prod_{l=1}^m\frac{(1-\beta)^{\rho_l}\alpha^{n{\rho_l}}}{(\alpha+l)^{n{\rho_l}}}b_l^{\rho_l}c_l^{\rho_l}\right)\right|\le\frac{2(1-\beta)\alpha^{n-1}}{(\alpha+v)^n}
\]
that is,
\begin{equation}
\left|\sum_{m=1}^v (-1)^{m+1}
\frac{(1-\beta)^m\alpha^{mn-1}\gamma_{m-1}}{2^m}\left(\sum_{j=1}^v
C_j\prod_{l=1}^m\frac{b_l^{\rho_l}c_l^{\rho_l}}{(\alpha+l)^{n{\rho_l}}}\right)\right|\le\frac{2(1-\beta)\alpha^{n-1}}{(\alpha+v)^n}\,.
\label{21}
\end{equation}
Now comparing ~(\ref{15}) and the term in the absolute value in
~(\ref{21}) (with $v=k-1$), we would conclude that the
inequalities ~(\ref{3}) hold if we are able to find two members
$h(z)=1+d_1z+d_2z^2+...$ and $H(z)=1+c_1z+c_2z^2+...$ of $P$ which
give rise to the constants $\gamma_m$ (as required by ~(\ref{6}))
and $c_l$. For $H\in P$, a natural choice of the Moebius function
is suitable. That is, $H(z)=(1+z)/(1-z)=1+2z+2z^2+...$. Thus we
have $c_l=2,\;\;\;l=1,\;2,\;...$. Using this in ~(\ref{21}) we get
\begin{equation}
\left|\sum_{m=1}^{k-1} (-1)^{m+1}
(1-\beta)^m\alpha^{mn-1}\gamma_{m-1}\left(\sum_{j=1}^{k-1}
C_j\prod_{l=1}^m\frac{b_l^{\rho_l}}{(\alpha+l)^{n{\rho_l}}}\right)\right|\le\frac{2(1-\beta)\alpha^{n-1}}{(\alpha+k-1)^n}\,.
\label{22}
\end{equation}
Comparing ~(\ref{15}) and the term in the absolute value in
~(\ref{22}) we find that
\[
(-1)^{m+1}\frac{\gamma_{m-1}}{\alpha}=\frac{\prod_{j=0}^{m-1}(1-j\alpha)}{m!\alpha^m}
\]
that is
\begin{equation}
\gamma_{m-1}=\frac{\prod_{j=1}^{m-1}(j\alpha-1)}{m!\alpha^{m-1}},\;\;\;\gamma_0=1\,.
\label{23}
\end{equation}
Now with $d_{\mu}$, $\mu=1,2,...,m-1$ defined by
\begin{equation}
\frac{1}{2^{m-1}}\left[1+\frac{1}{2}\sum_{\mu=1}^{m-1}\binom{m-1}{\mu}d_{\mu}\right]=\frac{\prod_{j=1}^{m-1}(j\alpha-1)}{m!\alpha^{m-1}},\,
\label{24}
\end{equation}
we need to find $h(z)_k$, corresponding to each $a_k$,
$k=2,3,4,...$., such that the coefficients $d_{\mu}$,
$\mu=1,2,...,m-1$ of each $h(z)_k$ satisfy ~(\ref{24}) and the
proof would be complete. Observe from ~(\ref{22}) that
$m=1,2,...,k-1$ as we begin to implement our scheme for each
$k=2,3,4,...$.

$k=2$: In this case $m=1$. Hence by ~(\ref{23}) $\gamma_0=1$ and
we can therefore define $d_{\mu}=0$ for all $\mu$ so that
$h(z)_2=1$.

$k=3$: Here we have $m=1,2$. Using ~(\ref{24}) we get
$d_1=-\tfrac{2}{\alpha}$, so that
\[
h(z)_3=\frac{\alpha-1}{\alpha}+\frac{1}{\alpha}\left(\frac{1-z}{1+z}\right)=1-\frac{2}{\alpha}z+...
\]

$k=4$: Now $m=1,2,3$. From ~(\ref{24}) we have
\[
\frac{1}{4}\left(1+\frac{1}{2}(2d_1+d_2)\right)=\frac{(\alpha-1)(2\alpha-1)}{3!\alpha^2}.
\]
Taking $d_1=0$, we get
\[
\frac{d_2}{2}=\frac{\alpha^2-6\alpha+2}{3\alpha^2}.
\]
(Note that $|d_2|\le 2$) and we have
\[
\begin{split}
&h(z)_4= \left\{
\begin{array}{ll}
\frac{2(\alpha^2+3\alpha-1)}{3\alpha^2}+\frac{\alpha^2-6\alpha+2}{3\alpha^2}\left(\frac{1+z^2}{1-z^2}\right)&\mbox{if $\alpha\ge \alpha_0$},\\
\frac{2(2\alpha^2-3\alpha+1)}{3\alpha^2}+\frac{6\alpha-\alpha^2-2}{3\alpha^2}\left(\frac{1-z^2}{1+z^2}\right)&\mbox{if
$1<\alpha\le \alpha_0$.}
\end{array}
\right.\\
&=1+\frac{2(\alpha^2-6\alpha+2)}{3\alpha^2}z^2+...
\end{split}
\]
where $\alpha_0>1$ is the solution of $\alpha^2-6\alpha+2=0$.

$k=5$: In this case $m=1,2,3,4$., and from ~(\ref{24}) we get
\[
\frac{1}{8}\left(1+\frac{1}{2}(3d_1+3d_2+d_3)\right)=\frac{(\alpha-1)(2\alpha-1)(3\alpha-1)}{4!\alpha^3}.
\]
Taking $d_1=d_2=0$, we get
\[
\frac{d_3}{2}=\frac{3\alpha^3-11\alpha^2+6\alpha-1}{3\alpha^3}.
\]
(Note also that $|d_3|\le 2$) and
\[
\begin{split}
&h(z)_5= \left\{
\begin{array}{ll}
\frac{11\alpha^2-6\alpha+1}{3\alpha^3}+\frac{3\alpha^3-11\alpha^2+6\alpha-1}{3\alpha^3}\left(\frac{1+z^3}{1-z^3}\right)&\mbox{if $\alpha\ge \alpha_0$},\\
\frac{6\alpha^3-11\alpha^2+6\alpha-1}{3\alpha^3}+\frac{1-6\alpha+11\alpha^2-3\alpha^3}{3\alpha^3}\left(\frac{1-z^3}{1+z^3}\right)&\mbox{if
$1<\alpha\le \alpha_0$.}
\end{array}
\right.\\
&=1+\frac{2(3\alpha^3-11\alpha^2+6\alpha-1)}{3\alpha^3}z^3+...
\end{split}
\]
In this case $\alpha_0>1$ is the solution of
$3\alpha^3-11\alpha^2+6\alpha-1=0$.

$k\ge 6$: In general, $m=1,2,3,...,k-1$. In ~(\ref{24}) we set
$d_1=\tfrac{-2}{m-1}$, $d_2=d_4=...=d_\xi=\sigma$ where $\xi$
equals $m-1$ if $m-1$ is {\em even} and $m-2$ otherwise, and
$d_3=d_5=...=d_\omega=0$ where $\omega$ equals $m-1$ if $m-1$ is
{\em odd} and $m-2$ otherwise. With this we get
\[
\frac{\sigma}{2}=\frac{2^{m-1}\prod_{j=1}^{m-1}\left(\frac{j\alpha-1}{j\alpha}\right)}{m\left(\binom{m-1}{2}+\binom{m-1}{4}+...+\binom{m-1}{\xi}\right)}.
\]
In this sense $|d_{\mu}|\le 2$ for all $\mu=1,2,...,m-1$. Setting
$m=k-1$, we now define $h(z)_k$, $k\ge 6$ as follows
\begin{multline*}
h(z)_k=\left(1-\frac{2}{k-2}-\frac{2^{k-2}\prod_{j=1}^{k-2}\left(\frac{j\alpha-1}{j\alpha}\right)}
{(k-1)\left(\binom{k-2}{2}+\binom{k-2}{4}+...+\binom{k-2}{\xi}\right)}\right)+\frac{2}{k-2}
(1-z)\\
+\frac{2^{k-2}\prod_{j=1}^{k-2}\left(\frac{j\alpha-1}{j\alpha}\right)}{(k-1)\left(\binom{k-2}{2}
+\binom{k-2}{4}+...+\binom{k-2}{\xi}\right)}\left(\frac{1+z^2}{1-z^2}\right).
\end{multline*}
In other words,
\begin{multline*}
h(z)_k=1-\frac{2}{k-2}z+\frac{2^{k-2}\prod_{j=1}^{k-2}\left(\frac{j\alpha-1}{j\alpha}\right)}
{(k-1)\left(\binom{k-2}{2}+\binom{k-2}{4}+...+\binom{k-2}{\xi}\right)}z^2\\
+\frac{2^{k-2}
\prod_{j=1}^{k-2}\left(\frac{j\alpha-1}{j\alpha}\right)}{(k-1)\left(\binom{k-2}{2}+\binom{k-2}{4}+...
+\binom{k-2}{\xi}\right)}z^4+...
\end{multline*}
That these functions $h(z)_k$ belong to $P$ follows from the fact
that the function $\lambda_1f_1+\lambda_2f_2+...+\lambda_mf_m$
belongs to $P$ if $f_1$, $f_2$, ..., $f_m$ belong, $\lambda_1$,
$\lambda_2$,...,$\lambda_m\ge 0$ and
$\lambda_1+\lambda_2+...+\lambda_m=1$. This completes the proof of
the theorem.

\begin{remark} We compute $h(z)_6$, ..., $h(z)_{10}$ for the purpose of illustration.
\begin{multline*}
h(z)_6=1-\frac{1}{2}z+\frac{4(\alpha-1)(2\alpha-1)...(4\alpha-1)}{105\alpha^4}z^2\\
\frac{4(\alpha-1)(2\alpha-1)...(4\alpha-1)}{105\alpha^4}z^4+...
\end{multline*}
\begin{multline*}
h(z)_7=1-\frac{2}{5}z+\frac{4(\alpha-1)(2\alpha-1)...(5\alpha-1)}{675\alpha^5}z^2\\
\frac{4(\alpha-1)(2\alpha-1)...(5\alpha-1)}{675\alpha^5}z^4+...
\end{multline*}
\begin{multline*}
h(z)_8=1-\frac{1}{3}z+\frac{8(\alpha-1)(2\alpha-1)...(6\alpha-1)}{10765\alpha^6}z^2\\
\frac{8(\alpha-1)(2\alpha-1)...(6\alpha-1)}{10765\alpha^6}z^4+...
\end{multline*}
\begin{multline*}
h(z)_9=1-\frac{2}{7}z+\frac{2(\alpha-1)(2\alpha-1)...(7\alpha-1)}{19845\alpha^7}z^2\\
\frac{2(\alpha-1)(2\alpha-1)...(7\alpha-1)}{19845\alpha^7}z^4+...
\end{multline*}
\begin{multline*}
h(z)_{10}=1-\frac{1}{4}z+\frac{4(\alpha-1)(2\alpha-1)...(8\alpha-1)}{360045\alpha^8}z^2\\
\frac{4(\alpha-1)(2\alpha-1)...(8\alpha-1)}{360045\alpha^8}z^4+...
\end{multline*}

With this work the coefficient problem of functions in the class
$T_n^\alpha(\beta)$ is settled for any $\alpha>1$. Of course the
case $\alpha=1$ is trivial as this simply gives
$|a_k|\le\tfrac{2(1-\beta)}{(k+1)^n}$, $k\ge 2$, as can be seen
easily from ~(\ref{12}). Thus the problem only remains open for
$(k-3)^{-1}\le\alpha<1$, $k\ge 5$. Finally, we note a humble
attempt at this problem made by the authors in \cite{OBFR}. Their
results depended wholly on the triangle inequality, and were not
sharp.
\end{remark}
 \medskip

{\it Acknowledgements.} Special thanks to the Abdus Salam
International Centre for Theoretical Physics, Trieste, Italy for
providing the reference paper \cite{NN} and particularly Dr.
Siraaj Ajadi, Department of Mathematics, Obafemi Awolowo
University, Ile-Ife, Nigeria for his {\em being there always}
during the author's postgraduate studies at Ife and Ilorin.

\vspace{10pt}

\hspace{-4mm}{\small{Received}}

% ADDRESS
\vspace{-12pt}
\ \hfill \
\begin{tabular}{c}
{\small\em  Department of Mathematics}\\
{\small\em  University of Ilorin}\\
{\small\em  Ilorin, Nigeria}\\
{\small\em E-mail: {\tt abuuabdilqayyuum@gmail.com}} \\
\end{tabular}

\end{document}